\documentclass[twoside,10pt]{article}
\usepackage{amssymb,amsmath,euscript}
\usepackage{graphicx,color,caption}
\usepackage{subfigure}
\usepackage{amsfonts,amssymb,amsthm,graphicx,booktabs,enumerate,mathrsfs,subfigure,multicol,mathdots}
\usepackage{epstopdf}
\usepackage[all]{xy}
\usepackage{cite}
\newtheorem{proposition}{Proposition}[section]
\newtheorem{theorem}[proposition]{Theorem}
\newtheorem{lemma}[proposition]{Lemma}

\newtheorem{definition}[proposition]{Definition}
\newtheorem{remark}[proposition]{Remark}
\newtheorem{example}[proposition]{Example}
\newtheorem{assumption}[proposition]{Assumption}

\RequirePackage{fix-cm}
\usepackage{graphicx}
\usepackage{multirow}
\usepackage{array}
\newcommand{\PreserveBackslash}[1]{\let\temp=\\#1\let\\=\temp}
\newcolumntype{C}[1]{>{\PreserveBackslash\centering}p{#1}}
\newcolumntype{R}[1]{>{\PreserveBackslash\raggedleft}p{#1}}
\newcolumntype{L}[1]{>{\PreserveBackslash\raggedright}p{#1}}
\usepackage{color}
\usepackage{subfigure,epstopdf}
\usepackage{amsmath,amsfonts,amssymb}
\usepackage{cases}
\usepackage{algorithm,algorithmic}
\usepackage{footmisc}
\usepackage{url}
 
\usepackage{breakurl}

\renewcommand{\S}{\mathbb{S}}
\newcommand{\T}{\mathbb{T}}
\newcommand{\R}{\mathbb{R}}
\newcommand{\J}{\mathcal{J}}
\newcommand{\N}{\mathcal{N}}

\title{\bf{A Lagrange-Newton Algorithm for Sparse Nonlinear Programming}\thanks{This research was partially supported by the National Natural Science Foundation of China (11771038, 11971052, 12011530155) and Beijing Natural Science Foundation (Z190002).}
\author{
Chen Zhao \thanks{Department of Mathematics, Beijing Jiaotong University, Beijing 100044, China; {\tt 14118409@bjtu.edu.cn}.} ,\hspace{1mm} 
Naihua Xiu \thanks{Department of Mathematics, Beijing Jiaotong University, Beijing 100044, China; {\tt nhxiu@bjtu.edu.cn}.} ,\hspace{1mm} 
Houduo Qi \thanks{ School of Mathematics, University of Southampton, Southampton SO17 1BJ, UK; {\tt hdqi@soton.ac.uk}.} ,\hspace{1mm} 
Ziyan Luo \thanks{Corresponding author; Department of Mathematics, Beijing Jiaotong University, Beijing 100044, China; {\tt zyluo@bjtu.edu.cn}.} 
}}

\begin{document}
\maketitle

\begin{abstract}
The sparse nonlinear programming (SNP) problem has wide applications in signal and image processing, machine learning, pattern recognition, finance and management, etc. However, the computational challenge posed by SNP has not yet been well resolved due to the nonconvex and discontinuous $\ell_0$-norm involved. In this paper, we resolve this numerical challenge by developing a fast Newton-type algorithm. As a theoretical cornerstone, we establish a first-order optimality condition for SNP based on the concept of strong $\beta$-Lagrangian stationarity via the Lagrangian function, and reformulate it as a system of nonlinear equations called the Lagrangian equations. The nonsingularity of the corresponding Jacobian is discussed, based on which the Lagrange-Newton algorithm (LNA) is then proposed. Under mild conditions, we establish the locally quadratic convergence and the iterative complexity estimation of LNA. To further demonstrate the efficiency and superiority of our proposed algorithm, we apply LNA to solve two specific application problems arising from compressed sensing and sparse high-order portfolio selection, in which significant benefits accrue from the restricted Newton step in LNA.

\end{abstract}

\section{Introduction}
In this paper, we are mainly concerned with the following sparse nonlinear programming (SNP) problem:
\begin{equation}\label{min}
\min\limits_{x\in\R^n} f(x),~~\text{s.t.}~h(x)=0,~x\in \S,
\end{equation}
where $f:\R^n\rightarrow\R$ and $h:=(h_1,\ldots,h_m)^{\top}:\R^n\rightarrow\R^m$  are twice continuously differentiable functions, $\S:=\{x\in \R^n:\|x\|_0\leq s\}$ is the sparse constraint set with integer $s\in (0,n)$ and $\|x\|_0$ is the $\ell_0$-norm of $x$ that counts the number of nonzero components of $x$. By denoting $\Omega:= \{x\in \R^n:h(x)=0\}$, the feasible set of problem \eqref{min} is abbreviated as $\Omega\cap\S$ and the optimal solution set can be written as $\arg\min_{x\in \Omega\cap \S}f(x)$. The SNP problem has wide applications ranging from linear and nonlinear compressed sensing \cite{Donoho2006Compressed,elad2010role} in signal processing, the sparse portfolio selection\cite{gao2013optimal,xu2015} in finance, to variable selection \cite{MisraInteractive,koh2007} and sparse principle component analysis \cite{ZouSparse,BeckThe} in high-dimensional statistical analysis and machine learning, etc. Unfortunately, due to the intrinsic combinatorial property in $\S$, the SNP problem is generally NP-hard, even for the simple convex quadratic objective function \cite{Nat1995}.

To well resolve the computational challenge resulting from $\S$, efforts have been made in two mainstreams in the literature. The first mainstream is ``{\it relaxation}" approach, with a rich variety of relaxation schemes distributed in
\cite{Chen2014Complexity,FanVariable,Gotoh2018}, just name a few.
The second one is the ``{\it greedy}" approach that tackles the involved $\ell_0$-norm directly, with a large number of algorithms tailored for SNP with the feasible set $\S$ merely (i.e., $\Omega=\R^n$), see, e.g., the first-order algorithms
\cite{tropp2007signal,needell2009cosamp,BLUMENSATH2009265},
and second-order algorithms with the Newton-type steps interpolated
\cite{Yuan2014Newton,YuanPAMI2017,JMLR:v18:14-415,10.1145/3097983.3098165}, etc. As the first-order information such as gradients are used in first-order greedy algorithms, linear rate convergence results are established as one can expect. While benefitting from the second-order information such as Hessian matrices, the aforementioned second-order greedy algorithms are witnessed in the numerical experiments with superior performance in terms of fast computation speed and high solution accuracy.
Besides the notable computational advantage that observed numerically, in a very recent work \cite{zhou2019global}, Zhou et al. propose a new algorithm called Newton Hard-Thresholding Pursuit (NHTP) with cheap Newton steps in a restricted fashion and rigorously establish the quadratic convergence rate.

In sharp contrast to the fruitful computational algorithms for nonlinear programming with the single sparse constraint set $\S$, a small portion of research on greedy algorithms is addressed for general SNP over $\S$ intersecting with some additional constraint set. The limited works are distributed in \cite{kyrillidis2012sparse,zhou2017IIHT,beck2015minimization,lu2015optimization,Lusparsepenalty}.
It is noteworthy that these algorithms are mostly gradient-based, and no quadratic convergence rate can be expected. To make up such a deficiency, the appealing theoretical and computational properties of NHTP \cite{zhou2019global} inspires us to develop a quadratic convergent Newton-type method for SNP when $\Omega$ is characterized by nonlinear equality constraints as presented in problem \eqref{min}.

The main contributions of this paper are summarized as below:
\begin{itemize}
\item[(i)] The strong $\beta$-Lagrangian stationarity is introduced to characterize the optimality condition for SNP, and an equivalent characterization of such a stationary point is built which is accessible to performing the Newton method.
\item[(ii)] The crucial Jacobian nonsingularity of the underlying system is well addressed under some mild assumptions, and the essential linear system in each iteration is reduced to be of size $(s+m)\times (s+m)$, a significant dimension reduction benefitting from the intrinsic sparsity.
\item[(iii)] The resulting Lagrange-Newton algorithm (LNA) is shown to possess the locally quadratic convergence, and to gain high efficiency  numerically for specific application problems including compressed sensing and sparse portfolio selection. 
\end{itemize}

The remainder of this paper is organized as follows. In Section 2, 
the optimality condition in terms of Lagrangian stationarity is established. In Section 3, an equivalent Lagrangian equation system is proposed
and the nonsingularity of its Jacobian is discussed.
The framework of LNA and its locally quadratic convergence are elaborated in Section 4.
Two well-known applications are analyzed in Section 5.
Extensive numerical experiments are conducted in Section 6. Conclusions are made in Section 7.

For convenience, the following notations will be used throughout the paper. For any given positive integer $n$, denote $[n]:=\{1,\ldots,n\}$.
For an index set $J \subseteq [n]$, let $|J|$ be the cardinality of $J$ that counts the number of elements in $J$, and $J^c:=[n]\setminus T$ be its complementary set. The collection of all index sets with cardinality $s$ in $[n]$ is defined by $\J_s:=\{J\subseteq [n]: |J| = s\}$. Given $x\in\R^n$,  
denote $\text{supp}(x):=\{i\in [n]:x_i\neq 0\}$ and $\J_s(x):=\{J\in\J_s: \text{supp}(x)\subseteq J\}$.
 We define $x_T\in\R^{|T|}$ as the subvector of $x\in\R^n$ indexed  by $T$. For the matrix $A\in\R^{m \times n}$, define $A_{I,J}$ as a submatrix whose rows and columns are respectively indexed by $I$ and $J$. In particular, we write $A_T$ as its submatrix consisting of columns indexed by $T$ and $A_{T,\cdot}$ as its submatrix consisting of rows indexed by $T$.
Given a twice continuously differentiable function $g$ with its gradient $\nabla g(x)$ and Hessian $\nabla^2 g(x)$ at $x$, denote $\nabla_{T}g(x):=(\nabla g(x))_{T}$
 and $\nabla^2_{I,J} g(x):=(\nabla^2 g(x))_{I,J}$.
The Euclidean norm of a vector $x$ is denoted by $\|x\|$, and the spectral norm of a matrix $A$ is denoted by $\|A\|$.

\section{Lagrangian Stationarity}

This section is devoted to the optimality conditions for \eqref{min} in terms of the Lagrangian stationarity, which will build up the theoretical fundamentals to our new proposed algorithm in the sequel.

Firstly, we consider the projection on sparse set $\S$. For any given nonempty closed set $Q\subseteq \R^n$ and any $x\in \R^n$, define the projection operator $\Pi_Q(x):=\arg\min_{y\in Q} \|x-y\|^2$. Recall from \cite{beck2015minimization,Pan2015} that the projection operator $\Pi_{\S}$ admits an explicit formula as follows: for any $z\in\R^n$, and for any  $\pi\in\Pi_{\S}(z)$, we have
\begin{equation}\label{explicit-proj-s}
\pi_{t_i}=
\left\{
\begin{array}{ll}
     z_{t_i},  &i\in [s],\\
      0,         &\text{otherwise},
    \end{array}
 \right.
\end{equation}
where $\{t_1,\ldots,t_n\}$ satisfies
$|z_{t_1}|\geq\ldots\geq |z_{t_n}|$. Utilizing the sparse projection $\Pi_{\S}$, Beck and Eldar  \cite[Theorem 2.2]{beck2013sparsity} introduced and characterized an optimality condition for SNP with single sparse contraint set $\S$: if $x^*$ is a global minimizer, then for any $L>L_f$, $x^*=\Pi_{\S}(x^*-\nabla f (x^*)/L)$, where $L_f$ is the Lipschitz constant of $\nabla f$. Later, Lu \cite{lu2015optimization} extended such a result to the case of $\Omega'\cap \S$  by means of 
$x^*=\Pi_{\Omega' \cap\S}(x^*- \nabla f (x^*)/L)$, where $\Omega'$ is a nonempty closed and convex set. Limitation follows when $\Pi_{\Omega' \cap\S}(\cdot)$ has no explicit expression. This motivates us to introduce the following Lagrangian stationarity.

\begin{definition}\label{sta}
 Given $x^*\in \R^n$ and $\beta>0$, $x^*$ is called a strong $\beta$-Lagrangian stationary point of problem \eqref{min} if there exists a Lagrangian multiplier $y^*\in \R^m$ such that
  \begin{equation}\label{strongbeta}
  \left\{
    \begin{array}{ll}
      x^*=\Pi_{\S}(x^*-\beta \nabla_xL (x^*,y^*)),\\
      h(x^*)=0,
    \end{array}
  \right.
  \end{equation}
where $L(x,y):=f(x)-\langle y, h(x)\rangle$ is the Lagrangian function associated with \eqref{min} for any $x\in \S$ and $y\in \R^m$.
\end{definition}

By employing  \eqref{explicit-proj-s} of $\Pi_{\S}$, one can easily rewrite system \eqref{strongbeta} as follows.

\begin{lemma}\label{exbeta} Given $x^*\in \R^n$ and $y^*\in \R^m$, denote $\Gamma^*:=supp(x^*)$ and $q^*=\nabla_xL(x^*,y^*)$.
Then $x^*$ is a strong $\beta$-Lagrangian stationary point of problem \eqref{min} with $y^*$ if and only if
\begin{numcases}{x^*\in \S,~h(x^*)=0,~  }
       \beta \|q^*_{(\Gamma^*)^c}\|_\infty <|x^*|_{(s)}  ~\&~ q^*_{\Gamma^*}=0,
  & \hbox{if $\|x^*\|_0 = s$,}\label{lemma2.2_21}\\
  q^*= 0,
  & \hbox{if $\|x^*\|_0 < s$,}\label{lemma2.2_22}
\end{numcases}
where $\text{supp}(x):=\{i\in [n]:x_i\neq 0\}$ and $x_{(s)}$ is the $s$th largest component of $x$.
\end{lemma}

The following assumption is introduced, followed by optimality analysis.

\begin{assumption}\label{assumption1} Given $x\in \Omega\cap\S$, $\text{rank}(\nabla_{\Gamma} h(x))=m$, where $\nabla h(x): =\left(\nabla h_1(x), \ldots, \nabla h_m(x)\right)^{\top}$, $\Gamma = \text{supp}(x)$.
\end{assumption}

\begin{theorem}[First-order necessary optimality condition]\label{optimalsta} Suppose that $x^*$ is a local minimizer of \eqref{min} and Assumption 1 holds at $x^*$. Then there exists a unique $y^*\in\R^m$ such that  $x^*$ is a strong $\beta$-Lagrangian stationary point of problem \eqref{min} for any $\beta\in (0, \hat{\beta})$, where
\begin{equation}\label{hatbeta}
\hat{\beta}:=\left\{ \begin{array}{ll}
\frac{|x^*|_{(s)}}{\left\|q^*_{(\Gamma^*)^c}\right\|_{\infty}},    &\mbox{~if~}\|x^*\|_0=s\mbox{~and~} q^*_{(\Gamma^*)^c}\neq 0,\\
+\infty, &\mbox{~otherwise.}
\end{array}\right.
\end{equation}
\end{theorem}
\noindent{\bf Proof.} Since $x^*$ is a local minimizer of (1.1), for any $J\in\J_s(x^*)$, $x^*$ is also a local  minimizer of
\begin{equation*}
(P_J)~~~\min\limits_{x\in\R^n} f(x),~~\text{s.t.}~h(x)=0,~x_{J^c}=0.
\end{equation*}
Assumption 1 implies that $\nabla h_i(x^*),~i\in [m]$, $e_j,~j\in J^c$ are linearly independent, where  $e_j\in\R^n$ is the $j$th column in the identity matrix. It means that linear independent constraint qualification (LICQ) holds at $x^*$ for ($P_J$). Thus for any given $J\in\J_s(x^*)$, there exists a unique $y^J\in\R^m$ and a unique $z^J\in\R^{n-s}$ such that
\begin{equation}\label{opconPJ}
\nabla f(x^*)=(\nabla h(x^*))^{\top}y^J+\sum\limits_{j\in J^c} z_j^J e_j.
\end{equation}

\noindent {Case I:} When $\|x^*\|_0=s$, we have $\J_s(x^*)=\{\Gamma^*\}$. Set $J = \Gamma^*$ in \eqref{opconPJ}, and let $y^*=y^{\Gamma^*}$.  
Direct calculations yield
$$q^*_{\Gamma^*}=\left(\nabla f(x^*)-(\nabla h(x^*))^{\top}y^*\right)_{\Gamma^*}=\left(\sum\limits_{i\in {(\Gamma^*)}^c} z_i^{\Gamma^*} e_i\right)_{\Gamma^*}=0 ~\mbox{~and~~} q^*_{(\Gamma^*)^c}=z^{\Gamma^*},$$
from which \eqref{lemma2.2_21} holds for all $\beta\in (0,\hat{\beta})$.

\noindent {Case II:} When $\|x^*\|_0<s$, $|\J_s(x^*)|={\mathcal{C}}_{n-\Gamma^*}^{s-\Gamma^*}$ is finite. For any $J\in\J_s(x^*)$, rewrite \eqref{opconPJ} into the block form
\begin{equation}\label{opcons}
\left\{
  \begin{array}{ll}
    \nabla_{\Gamma^*} f(x^*) =  (\nabla_{\Gamma^*}  h(x^*))^{\top}y^J, & \hbox{ } \\
    \nabla_{J\backslash\Gamma^*} f(x^*) = \nabla_{J\backslash\Gamma^*}  h(x^*))^{\top}y^J, & \hbox{ } \\
    \nabla_{J^c} f(x^*) = \nabla_{J^c}  h(x^*))^{\top}y^J + z^J. & \hbox{ }
  \end{array}
\right.
\end{equation}
By virtue of Assumption 1, the first equation in \eqref{opcons} indicates that $y^J$'s coincide for all $J\in\J_s(x^*)$, which we assign to $y^*$. It then leads to $q^*_{\Gamma^*}=0$.
As $\bigcup_{J\in\J_s(x^*)}(J\backslash\Gamma^*)= (\Gamma^*)^c$, it follows from the second equation in \eqref{opcons} that $q^*_{(\Gamma^*)^c} = 0$. Thus, $q^*=0$. This completes the proof by utilizing Lemma \ref{exbeta}. \qed

\begin{theorem}[First-order sufficient optimality condition]\label{sufficient-opt}
Let $f$ be a convex function and $h$ be an affine function. Given $\beta>0$, suppose that $x^*$ is a strong $\beta$-Lagrangian stationary point of \eqref{min} with the Lagrangian multiplier $y^*\in\R^m$. If $\|x^*\|_0=s$, then $x^*$ is a local minimizer of \eqref{min}; If $\|x^*\|_0<s$,  then $x^*$ is a global minimizer of \eqref{min}.
\end{theorem}
\noindent{\bf Proof.} Under the hypotheses on $f$ and $h$, the Lagrangian function $L(x,y)$ is convex with respect to $x$. It follows that
\begin{equation}
L(x,y^*)\geq L(x^*,y^*)+\langle q^*,x-x^* \rangle, ~\forall ~x\in\Omega\cap \S.
\end{equation}
Using the facts $x, x^*\in \Omega\cap \S $, we have $L(x,y^*)=f(x),~L(x^*,y^*)=f(x^*)$. Since $x^*$ is a strong $\beta$-Lagrangian stationary point with $y^*$, if $\|x^*\|_0<s$, $\langle q^*,x-x^* \rangle = 0$ from \eqref{lemma2.2_22}. 
Then for any $x\in \Omega\cap \S$, $f(x)\geq f(x^*)$. Thus $x^*$ is a global minimizer.
If $\|x^*\|_0=s$, there exists a sufficiently small $\delta>0$ such that for any $x\in\N(x^*,\delta)\cap(\Omega\cap \S)$, $x_{(\Gamma^*)^c} =0$ and hence $(x-x^*)_{(\Gamma^*)^c}=0$. By invoking \eqref{lemma2.2_21}, it yields that
$$\langle q^*,x-x^* \rangle= \langle q^*_{\Gamma^*},\left(x-x^*\right)_{\Gamma^*}\rangle + \langle q^*_{(\Gamma^*)^c},\left(x-x^*\right)_{(\Gamma^*)^c}\rangle = 0.$$ Thus for any $x\in\N(x^*,\delta)\cap(\Omega\cap \S)$, $f(x)\geq f(x^*)$, which implies that $x^*$ is a local minimizer of \eqref{min}.\qed

\vspace{2mm}

\begin{remark}Consider  problem \eqref{min}. (i) By virtue of Theorem \ref{optimalsta} and Lemma \ref{exbeta}, we summarize the relations among the strong $\beta$-Lagrangian stationarity (Strong-$\beta$-LS), B-KKT point and C-KKT point in \cite[Definition 3.1]{pan2017optimality}, S-stationarity (S-stat) and M-stationarity (M-stat) in \cite[Definition 4.1]{Michal2016constraint} as below.   
\begin{eqnarray*}\label{relation5}\text{Local~minimizer~} & \overset{{\footnotesize{\text{Assumption~\ref{assumption1}}}}}{\Longrightarrow} & \text{~Strong-}\beta\text{-LS~} \forall \beta\in (0,\hat{\beta}) \nonumber\\
& &~~~~~~~~~~~ \Updownarrow\\
\text{~M-stat}~\Longleftarrow~~\text{~S-stat}~~\Longleftrightarrow & ~~\text{~C-KKT~point}~~\Longleftarrow &\text{~B-KKT~point} \nonumber 
 \end{eqnarray*}
(ii) Given a feasible solution $x$, we can show that Assumption \ref{assumption1}, the restricted Robinson constraint qualification (R-RCQ) in
\cite[Definition 3.1]{pan2017restricted}, and the cardinality constraints linear independence constraint qualification (CC-LICQ) in \cite[Definition 3.11]{Michal2016constraint} are equivalent. They are stronger than the restricted linear independent constraint qualification (R-LICQ) in \cite[Definition 2.4]{pan2017optimality} when $\|x\|_0=s$.
Learning from \cite{pan2017optimality}, together with (i) of this remark, we obtain that R-LICQ could also ensure the existence of strong $\beta$-Lagrangian stationary point, but hardly could guarantee the uniqueness of the corresponding Lagrangian multiplier as stated in Theorem \ref{optimalsta}.
\end{remark}

\section{Lagrangian Equations and Jacobian Nonsingularity}
In this section, we will present an equivalent reformulation for the Lagrangian stationarity in terms of nonlinear equations, and discuss the Jacobian nonsingularity of the resulting equation system.

\subsection{Lagrangian Equations}
The optimality conditions in terms of the strong $\beta$-Lagrangian stationary point, as established in Theorems \ref{optimalsta} and \ref{sufficient-opt}, provide a way of solving \eqref{min}.
As one knows that $\Pi_\S(\cdot)$ is not differentiable, the main challenge is how to tackle such a non-differentiability. By exploiting the special structure possessed by the projection operator $\Pi_{\S}$, we propose a differentiable reformulation of the definitional expression of the strong $\beta$-Lagrangian stationary point \eqref{strongbeta}, using a finite sequence of Lagrangian equations.

\begin{definition} Given $x\in\S$, $y\in\R^m$ and $\beta>0$, denote $u:= x-\beta \nabla_xL(x,y)$. Define the collection of sparse projection index sets of $u$ by
\begin{equation}\label{Txyb}
\T(x,y;\beta):=\{T\in \J_s:  |u_i|\geq |u_j|,~\forall i \in T, \forall j \in T^c\}.
\end{equation}
For any given $T \in \T (x,y;\beta)$,
define the corresponding Lagrangian equation as
\begin{equation}\label{Fxy=0}
F(x,y;T):=\left[  \begin{matrix}
     (\nabla_xL(x,y))_{T} \\
     x_{T^c} \\
     -h(x)
\end{matrix}   \right]
=0.
\end{equation}
\end{definition}

\vspace{2mm}
As one can see, the function $F(x,y;T)$ in \eqref{Fxy=0} is differentiable with respect to $x$ and $y$ once $T$ is selected. Moreover, we have the following equivalent relationship between \eqref{Fxy=0} and \eqref{strongbeta}.

\begin{theorem}\label{eqeq} Given $x^*\in \S$, $y^*\in\R^m$ and $\beta>0$, $x^*$ is a strong $\beta$-Lagrangian stationary point of \eqref{min}  with the Lagrangian multiplier $y^*$ if and only if for any $T\in \T(x^*,y^*;\beta)$, $F(x^*,y^*;T)=0$.
Meanwhile,
 $\T(x^*, y^*; \beta) = \J_s(x^*).$
\end{theorem}

\noindent{\bf Proof.} By invoking the proof Lemma 4 in \cite{zhou2019global}, we can obtain the equivalent relationship directly. Now, we prove the rest part of the  theorem. If follows from Theorem \ref{eqeq} that for any $T\in \T(x^*, y^*; \beta)$, $F(x^*, y^*; T) = 0$. Thus, $x^*_{T^c} =0$ and hence $\text{supp}(x^*)\subseteq T$. It then yields that $T\in \J_s(x^*)$. The arbitrariness of $T$ leads to the inclusion $\T(x^*, y^*; \beta)\subseteq \J_s(x^*)$. It now suffices to show $\J_s(x^*)\subseteq \T(x^*, y^*; \beta)$. If $\|x^*\|_0 <s$, then $\nabla_x L(x^*, y^*) = 0$ by invoking \eqref{lemma2.2_22}, and hence $u^* = x^* -\beta \nabla_x L(x^*, y^*) =x^*$. For any $T\in \J_s(x^*)$, it is easy to verify that for any $i\in T$ and any $j\in T^c$, $|u^*_i|=|x^*_i|\geq |x^*_j| = |u^*_j|$, which indicates that $T\in \T(x^*, y^*;\beta)$; If $\|x^*\|_0 =s$, then $|x^*|_{(s)} >0$ and $\J_s(x^*) =\{\Gamma^*\}$. By virtue of  \eqref{lemma2.2_21}, we have that for any $i\in \Gamma^*$ and any $j\notin \Gamma^*$,
$|u^*_i| = |x^*_i| \geq |x^*|_{(s)}> \beta |(\nabla_x L(x^*, y^*))_j| = |u^*_j|$. Thus, $\Gamma^*\in \T(x^*, y^*;\beta)$. In a word, in both cases, we can conclude $\J_s(x^*)\subseteq \T(x^*, y^*; \beta)$. This completes the proof. \qed

\subsection{Jacobian Nonsingularity}
In this subsection, let $x^*$ be a strong $\beta$-Lagrangian stationary point with $y^*$. To handle the Lagrangian equation \eqref{Fxy=0} for a given index set $T\in\T(x^*,y^*;\beta)$, it is crucial to discuss the nonsingularity of the Jacobian of $F(x,y;T)$ with respect to $(x,y)$ in a neighborhood of $(x^*,y^*)$, namely,
\begin{equation}\label{F'xyT}
\nabla_{(x,y)}F(x,y;T)=
\left[\begin{matrix}
(\nabla^2_{xx} L(x,y) )_{T,\cdot}    & -(\nabla_{T}  h(x))^{\top}\\
I_{T^c,\cdot}& 0\\
-\nabla h(x)         &0
\end{matrix}\right]\in\R^{(n+m)\times (n+m)},
\end{equation}
where $\nabla^2_{xx} L(x,y) =\nabla^2 f(x)-\sum\limits_{i=1}^{m}y_i \nabla^2 h_i(x)$ is Hessian matrix of $L(x,y)$ with respect to $x$,  and $I$ is the identity matrix.
It is worth mentioning that since $T$ is related to $(x,y)$, the conventional Jacobian of $F(x,y;T)$ may differ with \eqref{F'xyT} if we treat $T$ as a function of $(x,y)$. However, as $T$ may vary as $(x,y)$ changes, we will update such an index set $T$ in our proposed iterative algorithm adaptively. 
Two additional assumptions are stated as below.
\vskip 1mm

\noindent{\bf Assumption 1$'$} rank$(\nabla_T h(x^*) )= m$, for any $T\in \J_s(x^*)$.

\vskip 1mm

\begin{assumption}\label{assumption3} (Second-order optimality condition) For any $T\in \J_s(x^*)$, $(\nabla^2_{xx} L(x^*,y^*) )_{T,T}$ is positive definite restricted to the null space of $\nabla_T h(x^*)$, i.e.,
$$ d^\top\!\left(\nabla^2_{xx} L(x^*,y^*)\right)_{T,T} \!d>0,\forall~0\neq d\in {\mathbb{N}}(\nabla_T h(x^*))\!:=\!\{d\in\R^s\!:\!\nabla_T h(x^*)d=0\}.$$
\end{assumption}
Given any $T\in \J_s$ and $(x,y)\in \S\times \R^m$, elementary row operations yield the equivalence between the nonsingularity of $\nabla_{(x,y)}F(x,y;T)$ and that of
\begin{equation}\label{G-for}
G(x,y;T):= \left[\begin{matrix}
(\nabla^2_{xx} L(x,y) )_{T,T}     & -(\nabla_{T}  h(x))^\top\\
-\nabla_{T}  h(x) &0
\end{matrix}\right]\in \R^{(s+m)\times (s+m)}.
\end{equation}
Thus, we call $G(x,y;T)$ the {\it reduced Jacobian} of $F(x,y;T)$. Furthermore, when Assumptions 1$'$  and \ref{assumption3} hold, we can directly obtain the desired nonsingularity at $(x^*, y^*)$ as stated in the following theorem.

\begin{theorem}\label{A2B2} Let $x^*$ be a strong $\beta$-Lagrangian stationary point with $y^*$. If Assumptions 1$'$ and \ref{assumption3} hold,
then $\nabla_{(x,y)}F(x^*,y^*;T)$ is nonsingular for each index set $T\in\T(x^*,y^*;\beta)$.
\end{theorem}

\vskip 2mm

The rest of this subsection is devoted to the nonsingularity of $G(x,y;T)$ when $(x,y)$ are sufficiently close to $(x^*,y^*)$, by employing the achieved nonsingularity of $G(x^*,y^*;T)$ and the following assumption.
\begin{assumption}\label{assumption4} $\nabla^2 f$ and $\nabla^2 h_i ~(i\in[m])$ are Lipschitz continuous near $x^*$.
\end{assumption}

The locally Lipschitz continuity in Assumption \ref{assumption4} allows us to find positive constants $\delta_0^*$, $L_1$, $L_2$ such that  for any $z,\hat{z}\in {\N}(z^*,\delta^*_0)$ with $z^*:=(x^*;y^*)$, we have
\begin{eqnarray}\label{A2}
\|\nabla_x L(x,y)-\nabla_x L(\hat x,\hat y)\|\leq L_1 \|z-\hat z\|,\\\label{A5}
\|\nabla^2L(x,y)-\nabla^2 L(\hat{x},\hat{y})\|\leq L_2  \left\| z-\hat z  \right\|.
\end{eqnarray}
Let $x^*\neq 0$ (since the trivial case $x^*=0$ is not desired in practice) be a strong $\beta$-Lagrangian stationary point with $y^*$. We can define
$$\delta^*_1:=\frac{\min\limits_{i\in\Gamma^*}|x_i^*|-\beta  \max\limits_{i\in(\Gamma^*)^c}|q_i^*|}{\sqrt{2}(1+\beta L_1)},$$
where $\Gamma^*$, $q^*$ are defined in Lemma \ref{exbeta}.
By employing Lemma \ref{exbeta}, one can easily verify that $\delta^*_1>0$ since $x^*\neq 0$. Denote
\begin{equation}\label{delta*}
\delta^* := \min\{\delta^*_0, \delta^*_1\},
\end{equation}
\begin{equation}\label{NS}
\N_{\S}(z^*;\delta^*):=\{z\in\R^{n+m}:x\in \S,~\|z-z^*\|<\delta^*\}.
\end{equation}

\begin{lemma}\label{lemmaa}
Let $x^*$ be a strong $\beta$-Lagrangian stationary point with $y^*$. Denote $z^* := (x^*;y^*)$. If Assumption \ref{assumption4} holds, then for any $z := (x; y)\in\N_{\S}(z^*;\delta^*)$, we have
\begin{equation}\label{set-inclusion}
\T(x,y;\beta)\subseteq\T(x^*,y^*;\beta)~~\text{and}~~ \Gamma^*\subseteq \text{supp}(x)\cap T, \forall ~T\in \T(x,y;\beta).
\end{equation}
Particularly, if $\|x^*\|_0=s$, then $\{\text{supp}(x)\}=\T(x,y;\beta)=\T(x^*,y^*;\beta)=\{\Gamma^*\}$.

\end{lemma}
\noindent{\bf Proof.} Since $x^*$ is a strong $\beta$-Lagrangian stationary point with $y^*$, we have $F(x^*,y^*;T)=0$, $\forall~T\in \T(x^*,y^*;\beta)= \J_s(x^*)$ from Theorem \ref{eqeq}.
Consider any given $z=(x;y)\in \N_{\S}(z^*;\delta^*)$, denote $\Gamma = \text{supp}(x)$ and $q = \nabla_x L(x, y)$. For any $i\in\Gamma^*$ and any $j\in (\Gamma^*)^c$, we have
\begin{eqnarray*}
&&|x_i-\beta q_i|-|x_j-\beta q_j| \\
&\geq & |x_i^*|-|x_i-x_i^*|-|x_j-x_j^*|-\beta |q_i-q_i^*|-\beta |q_j-q_j^*|-\beta |q^*_j|\\
&\overset{\eqref{A2}}{\geq} &  \min\limits_{t\in\Gamma^*}|x_t^*|-\sqrt{2}\|z-z^*\|-\sqrt{2}\beta L_1\|z-z^*\|- \beta  \max\limits_{t\in (\Gamma^*)^c}|q_t^*|\\
& \geq &  \min\limits_{t\in\Gamma^*}|x_t^*|-\beta  \max\limits_{t\in (\Gamma^*)^c}|q_t^*|-\sqrt{2}(1+\beta L_1)\delta^*\\
& \geq & 0.
\end{eqnarray*}
This indicates that $i\in T$ and hence
\begin{equation}\label{gammastar}
\Gamma^*\subseteq T, ~~\forall ~ T\in \T(x,y;\beta).
\end{equation}
Furthermore, we have
\begin{equation}\label{jss}
\T(x,y;\beta) \subseteq \J_s(x^*) = \T(x^*,y^*;\beta).
\end{equation}
Next we claim that $\Gamma^*$ is also a subset of $\Gamma$. If not, there exists an index $i_0 \in \Gamma^*\setminus \Gamma$. Then
$$\|z-z^*\|\geq\|x-x^*\|\geq|(x-x^*)_{i_0}|= |x^*_{i_0}|\geq \min\limits_{t\in\Gamma^*}|x_t^*|\geq \delta^*,$$
which is a contradiction to $z\in \N_{\S}(z^*;\delta^*)$. Thus,
\begin{equation}\label{gamma}
\Gamma^*\subseteq \Gamma.
\end{equation}
Summarizing \eqref{gammastar}, \eqref{jss} and \eqref{gamma}, we get \eqref{set-inclusion}. Particularly, if $\|x^*\|_0 = s$, then $|\Gamma^*|=s$ and $\J_s(x^*)=\{\Gamma^*\}$. Utilizing \eqref{gammastar}, \eqref{jss} and \eqref{gamma} again, together with the fact $|\Gamma| \leq s$, we immediately get the rest of the desired assertion. \qed

\vspace{2mm}
Finally, the desired nonsingularity in a given neighborhood is stated.

\begin{theorem}\label{new} Let $x^*$ be a strong $\beta$-Lagrangian stationary point with $y^*$. If Assumptions 1$'$, \ref{assumption3} and \ref{assumption4} hold, then there exist constants $\tilde{\delta}^*\in (0,\delta^*]$ and $M^*\in (0, +\infty)$ such that for any $z:=(x;y) \in \N_{\S}(z^*;\tilde{\delta}^*)$ with $z^*:=(x^*;y^*)$, the reduced Jacobian matrix $G(x,y;T)$ is nonsingular and
\begin{equation}\label{inv-bound}
\|G^{-1}(x,y;T)\|\leq M^*,~~\forall T\in \T(x,y;\beta).
\end{equation}
\end{theorem}
\noindent{\bf Proof.} By invoking Theorems \ref{eqeq} and \ref{A2B2}, we can get the nonsingularity of $G(x^*, y^*; T)$ for all $T\in \T(x^*, y^*;\beta) = \J_s(x^*)$. For any $z \in \N_{\S}(z^*;\delta^*)$, the inclusion $\T(x,y;\beta)\subseteq\T(x^*,y^*;\beta)$ from Lemma \ref{lemmaa} immediately yields the nonsingularity of $G(x^*, y^*; T)$ for each $T\in \T(x, y;\beta)$. Furthermore, it follows from \eqref{A5} that for any $T\in \T(x, y;\beta)$,
\begin{equation}\label{GLip}
\|G(x,y;T)-G(x^*,y^*;T)\| \leq \|\nabla^2 L(x,y)-\nabla^2 L(x^*,y^*)\|\leq L_2 \|z-z^*\|,
\end{equation}
which indicates that $G(\cdot,\cdot~;T)$ is Lipschitz continuous near $z^*$ for any given $T\in \T(x, y;\beta)$. Thus, there exists $\delta_T>0$ and $M_T>0$ such that for any $z\in \N(z^*;\delta_T)$, $G^{-1}(x,y;T)$ exists and $\|G^{-1}(x,y;T)\|\leq M_T$. Set
\begin{equation}\label{delta-M}
\tilde{\delta}^* :=\min\{\delta^*, \{\delta_T\}_{T\in \J_s(x^*)}\}, \text{~and~} M^*:= \max\limits_{T\in \J_s(x^*)}\{M_T\}.
\end{equation}
It follows readily that for any $z\in \N(z^*;\tilde{\delta}^*)$, $G(x,y;T)$ is nonsingular and $\|G^{-1}(x,y;T)\|\leq M^*$, for all $T\in \T(x,y;\beta)$. \qed

\section{The Lagrange-Newton Algorithm}
In this section, we propose a Newton Algorithm for solving Lagrangian equation \eqref{Fxy=0} of problem \eqref{min} which is named as Lagrange-Newton Algorithm (LNA), and analyze the convergence rate of the algorithm.
\subsection{LNA Framework}

By employing the relationship between the strong $\beta$-Lagrangian stationary point and the Lagrangian equations as stated in Theorem \ref{eqeq},  the basic idea behind our algorithm is: solve the Lagrangian equation $F(x,y;$ $T) = 0$ iteratively by using the Newton method, and update the involved index set $T$ accordingly from $\T(x, y; \beta)$ by definition in each iteration. Details on the algorithm are as below.

Given $\beta>0$, let $(x^k,y^k)\in \S\times \R^m$ be the current iteration. 
\vskip 1mm

\noindent {\tt Index Set Selection:} Choose one index set $T_k$ from $\T(x^k,y^k;\beta)$ defined as in \eqref{Txyb}. This can be safely accomplished by picking the indices of the first $s$ largest elements (in magnitude) in $x^k- \beta \nabla_x L(x^k,y^k)$.
\vskip 1mm

\noindent{\tt The Newton Step:}
The classical Newton equation is 
\begin{equation}\label{newtonform}
\nabla_{(x,y)}F(x^k,y^k;T_k)(x^{k+1}-x^k; y^{k+1}-y^k)=-F(x^k,y^k;T_k).
\end{equation}
After simple calculations, \eqref{newtonform} can be rewritten as
\begin{equation}\label{finalnewton}
\left\{
\begin{array}{ll}
x^{k+1}_{T^c_k}=0;\\
G(x^k,y^k;T_k)
\left[\begin{matrix}
x^{k+1}_{T_k}\\
y^{k+1}
\end{matrix}\right]=
\left[\begin{matrix}
- \nabla_{T_k} f(x^k)+(\nabla^2_{xx} L(x^k,y^k) )_{T_k,.} x^k    \\
h(x^k)-\nabla h(x^k)x^k
\end{matrix}\right],
\end{array}
\right.
\end{equation} from which a significant dimension reduction is attained, from $(n+m)\times(n+m)$ to $(s+m)\times (s+m)$.
Under the conditions presented in Theorem \ref{new}, $G(x^k,y^k;T_k)$ 
is nonsingular, and hence the next iteration can be obtained from the unique solution of \eqref{finalnewton}, which can be solved in a direct way if $s+m$ is small or by employing the conjugate gradient (CG) method when $s+m$ is relatively large. As indicated, the low computational cost of the Newton step is greatly benificial from the intrinsic sparisity, especially when $s\ll n$.
\vskip 1mm

\noindent {\tt Stopping Criterion:} Given the current iteration triplet $(x^k,y^k,T_k)$, to measure how far $x^k$ is from being a strong $\beta$-Lagrangian stationary point, the following quantity is adopted
\begin{equation} \label{Tolerance-Function}
 \eta_{\beta}(x^k\!,y^k\!; T_k) \!:= \!\| F(x^k\!,y^k\!; T_k) \| + \max_{i \in T_k^c} \!\left\{ \!\max \!\Big( |(\nabla_x L(x^k,y^k) )_i| \!- \!|x^k|_{(s)}/\beta, \ 0  \Big) \right\}.
\end{equation}
The first term on the right-hand side of \eqref{Tolerance-Function} is to measure the residual of the Lagrangian equation system, and the second term is to testify $$\left\|\left(\nabla_x L(x^k, y^k)\right)_{ T_k^c}\right\|_{\infty} \leq\frac{1}{\beta}|x^k|_{(s)},$$
an inequality comes from \eqref{Txyb} and \eqref{Fxy=0} after simple manipulations. The stopping criterion is then designed in terms of $\eta_{\beta}(x^k\!,y^k\!; T_k)$.

\noindent The algorithmic framework is now summarized as follows.

\begin{algorithm*}
\caption{Lagrange-Newton Algorithm (LNA) for \eqref{min}}
{\bf Step 0.} {\tt (Initialization)} Give $\beta>0$ and $\epsilon>0$, choose $(x^0,y^0)\in \S\times \R^m$ and set $k=0$.

{\bf Step 1.} {\tt (Index Set Selection)} Choose $T_k\in\T(x^k,y^k;\beta)$ by \eqref{Txyb}.

 {\bf Step 2.}  {\tt (Stopping Criterion)} If $\eta_{\beta}(x^k,y^k;\; T_k)\leq \epsilon$,
then stop. Otherwise, go to {\bf Step 3.}

{\bf Step 3.} {\tt (The Newton Step)} Update  $(x^{k+1},y^{k+1})$ by \eqref{finalnewton}, set $k = k+1$ and go to {\bf Step 1.}

\end{algorithm*}

\subsection{Locally Quadratic Convergence}
The locally quadratic convergence of LNA is shown to be inherited from the classic Newton method, armed with the essential invariance property of index sets as stated in Lemma \ref{lemmaa}. Specifically, we have
\begin{theorem}\label{Thm-Qudratic-Convergence} Given $\beta>0$, suppose $x^*$ is a strong $\beta$-Lagrangian stationary point of \eqref{min} with $y^*$. If Assumptions 1$'$, \ref{assumption3} and \ref{assumption4} hold. Let $\tilde{\delta}^*$ and $M^*$ be defined as in \eqref{delta-M}, and $\T(x^*,y^*; \beta)$ be defined as in \eqref{Txyb}, respectively. Denote $z^*:=(x^*;y^*)$.
Suppose that the initial point $z^0:=(x^0; y^0)$ of LNA satisfies $z^0\in \N_{\S}(z^*,\delta)$
with $\delta = \min\{\tilde{\delta}^*,\frac{1}{M^*L_2}\}$. Then the sequence $\{z^k:=(x^k;y^k)\}$ generated by LNA is well-defined and for any $k\geq 0$,
\begin{itemize}
\item [(i)]$\lim\limits_{k\rightarrow \infty} z^k=z^*$ with quadratic convergence rate, namely
$$\|z^{k+1}-z^*\|\leq \dfrac{M^*L_2}{2}\left \| z^k-z^* \right\|^2.$$
\item [(ii)] $\lim\limits_{k\rightarrow \infty} F(z^k;T_{k})=0$ with quadratic convergence rate, namely
    $$\|F(z^{k+1};T_{k+1})\|\leq  \frac{M^*L_2\sqrt{L_1^2+1}}{\lambda_H}\|F(z^k;T_{k})\|^2, $$
where $\lambda_H :=\min\limits_{T_k\in \J_s(x^*)}\lambda_{\min}\left(\nabla_z F(z^*;T_k)^\top \nabla_z F(z^*;T_k)\right) $.
\item[(iii)] $\lim\limits_{k\rightarrow \infty} \eta_{\beta}(z^k;T_{k})=0$ with $\eta_{\beta}(z^{k+1};T_{k+1}) \leq M^*L_2\sqrt{L_1^2+1}\|z^k-z^*\|^2$ and LNA will terminate when $$ k\geq \left\lceil \frac{\log_2\left(4\delta^2 M^*L_2\sqrt{L_1^2+1}/\epsilon\right)}{2}\right\rceil,$$
where $\lceil t\rceil$ denotes the smallest integer no less than $t$.
\end{itemize}
\end{theorem}

\noindent{\bf Proof.} By employing Theorem \ref{new}, we know that the sequence $\{z^k\}$ generated by LNA is well-defined from the nonsingularity of $G(x^k, y^k;T_k)$ for all $k\geq 0$.

(i) Choose $T_0\in\T(x^0,y^0;\beta)$. Lemma \ref{lemmaa}, together with Theorem \ref{eqeq}, yields $x^*_{T_0^c}=0$. Meanwhile, following from \eqref{finalnewton} in Algorithm 1, we also have $x^1_{T_0^c}=0$.
With some routine work, one can further obtain 
{ \begin{eqnarray}\label{add-change}
&&\|z^1-z^*\|\nonumber\\
&\overset{(\ref{finalnewton})}{=}&
\left\|\left[\begin{matrix}
x^1_{ T_{0} }\\
y^1
\end{matrix}\right]-\left[\begin{matrix}
x^*_{ T_{0} }\\
y^*
\end{matrix}\right]\right\| \overset{(\ref{inv-bound})}{\leq} M^*\left\|G(x^0,y^0;T_0)\left(\left[\begin{matrix}
x^1_{ T_{0} }\\
y^1
\end{matrix}\right]-\left[\begin{matrix}
x^*_{ T_{0} }\\
y^*
\end{matrix}\right]\right)\right\|\nonumber\\
&\leq &  \frac{M^*L_2}{2} \left \| z^0-z^* \right\|^2 \leq \frac{\left \| z^0-z^* \right\|}{2}  < \frac{\delta}{2}.
\end{eqnarray}}
Similar reasons allow us to sequentially get
\begin{equation}\label{zk-z*}
\|z^{k+1}-z^*\|\leq  \dfrac{M^*L_2}{2}\left \| z^k-z^* \right\|^2\text{~and~}z^k\in \N_{\S}(z^*, \frac{\delta}{2^k}).
\end{equation}
Hence $\lim\limits_{k\rightarrow \infty} z^k=z^*$ with quadratic convergence rate.

(ii) Similar to the case $k=1$, Lemma \ref{lemmaa} and Theorem \ref{eqeq} also yield $x^*_{T_{k+1}^c}=0$ for any $k\geq 1$ from \eqref{zk-z*}.
After basic manipulations, we have
\begin{equation}\label{power4} \|F_{\beta}(z^{k+1};T_{k+1})\|^2
 \overset{\eqref{add-change}}{\leq} (L_1^2+1)\left( \dfrac{M^*L_2}{2}\right)^2\|z^{k}-z^*\|^4.
\end{equation}
Additionally, the index set property in Lemma \ref{lemmaa} ensures $$T_k\in \T(z^k;\beta)\subseteq \T(z^*;\beta)=\J_s(x^*), ~\forall~ k\geq 0$$ which further leads to $F(z^*;T_k)=0$ from Theorem \ref{eqeq}, and the nonsingularity of $\nabla_z F(z^*;T_k)$ from Theorem \ref{A2B2}.
Thus, $\lambda_H>0$ and
\begin{eqnarray}\label{F-inv}
\|z^k-z^*\|^2
& \leq & \frac{1}{\lambda_H} \|\nabla_z F(z^*;T_k)(z^k-z^*)\|^2 \nonumber\\
& \leq & \frac{2}{\lambda_H} \|\nabla_z F(z^*;T_k)(z^k-z^*)+o(\|z^k-z^*\|)\|^2 \nonumber\\
& =    & \frac{2}{\lambda_H} \|F(z^k;T_k)\|^2,
\end{eqnarray}
where the last equality is from $F(z^*;T_k)=0$. Combining with \eqref{power4}, we have
$$\|F(z^{k+1};T_{k+1})\| \leq \frac{M^*L_2\sqrt{L_1^2+1}}{2}\|z^k-z^*\|^2 \leq \frac{M^*L_2\sqrt{L_1^2+1}}{\lambda_H}\|F(z^k;T_k)\|^2.$$

(iii) Denote $\zeta_k: = \max\limits_{j\in T^c_k} \{\max\{|(\nabla_x L(x^k,y^k))_j|- |x^k|_{(s)}/\beta, 0\},$ for any given $k\geq 1$. We claim that
\begin{equation}\label{main-claim}
\zeta_k\leq L_1 \|z^k-z^*\|,~~\forall k\geq 1.
\end{equation}
For each given $k\geq 1$, we consider the following two cases.

\noindent Case I: If $\|x^*\|_0<s$, then $\nabla_x L(x^*,y^*) = 0$. It follows from the definition of $\zeta_k$ that
$$ \zeta_k \leq \|(\nabla_x L(x^k,y^k))_{T_k^c}\| \leq \|\nabla_x L(x^k,y^k) - \nabla_x L(x^*,y^*)\| \leq L_1 \|z^k-z^*\|.$$

\noindent Case II: If $\|x^*\|_0=s$, we have $\T(z^k;\beta) = \{\Gamma_k\} = \{\Gamma^*\}$ from Lemma \ref{lemmaa}. Thus, for any $k\geq 1$, $T_k = \Gamma^*$ and $(\nabla_x L(x^*,y^*))_{T_k} =0$.
Besides, since $|x^k|_{(s)}>0$, there exists $i_k\in T_k $ such that $|x^k_{i_k}| = |x^k|_{(s)}$. For any $k\geq 1$, it follows from the definition of $\T(x^k,y^k;\beta)$ that for any $j\in T_k^c = (\Gamma^*)^c$,
\begin{eqnarray}\label{ineq6}
\beta |(\nabla_x L(x^k,y^k))_j|
&  =  & |x_j^k - \beta(\nabla_x L(x^k,y^k))_j| \nonumber\\
&\leq & \min\limits_{i\in T_k}\{ |x_i^k - \beta(\nabla_x L(x^k,y^k))_i|\} \nonumber\\
&\leq & |x_{i_k}^k - \beta(\nabla_x L(x^k,y^k))_{i_k}|\nonumber\\
&  \leq  & |x^k|_{(s)} + \beta\|(\nabla_x L(x^k,y^k))_{T_k}\|.
\end{eqnarray}
This implies that $ \max\{|(\nabla_x L(x^k,y^k))_j|- |x^k|_{(s)}/\beta, 0\} \leq \|(\nabla_x L(x^k,y^k))_{T_k}\|, ~~\forall j\in T_k^c,$
which means $\zeta_k \leq \|(\nabla_x L(x^k,y^k))_{T_k}\|$. Together with $(\nabla_x L(x^*,y^*))_{T_k} =0$, we have
$$\zeta_k \leq \|(\nabla_x L(x^k,y^k))_{T_k}-(\nabla_x L(x^*,y^*))_{T_k}\|  \leq L_1 \|z^k-z^*\|,~~\forall k\geq 1.$$
This shows the claim in \eqref{main-claim}. Combining with \eqref{power4} and \eqref{add-change}, we further get that for $k\geq 1$,
\begin{eqnarray}
\eta_{\beta}(x^k,y^k;T_k)
&  =  & \|F(z^k;T_k)\|+\zeta_k \nonumber\\
&\leq & \frac{M^*L_2\sqrt{L_1^2+1}}{2}\|z^{k-1}-z^*\|^2 +\frac{M^*L_2L_1}{2}\|z^{k-1}-z^*\|^2\nonumber\\
&\leq &  M^*L_2\sqrt{L_1^2+1}\|z^{k-1}-z^*\|^2.
\end{eqnarray}
In addition, by virtue of \eqref{zk-z*}, we obtain $\eta_{\beta}(x^k,y^k;T_k) \leq \frac{\delta^2 M^*L_2\sqrt{L_1^2+1}}{2^{2k-2}}.$
To meet the stopping criterion $\eta_{\beta}(x^k,y^k;T_k)\leq \epsilon$ in LNA, it suffices to have $\frac{\delta^2 M^*L_2\sqrt{L_1^2+1}}{2^{2k-2}}\leq  \epsilon$, which leads to the bound of $k$ as desired.
This completes the proof. \qed

\vskip 3mm

\section{Applications}
Two selected SNP problems arising from some important applications are considered to demonstrate the effectiveness of our proposed Lagrange-Newton algorithm.

\subsection{Compressed Sensing}
Compressed sensing (CS) \cite{Michael2010Sparse} has been widely applied  in
signal and image processing \cite{elad2010role}, machine learning \cite{zhou2011manifold}, statistics
\cite{negahban2012unified}, etc.
A more general framework is considered, where some noise-free observations are allowed and added as hard constraints into the standard CS model,  taking the form of
\begin{equation}\label{CS12}
\min ~\frac{1}{2}\|Ax-b\|^2,~~\text{s.t.}~Cx=d,~x\in \S,
\end{equation}
where $A\in\R^{(p-m)\times n}$, $C\in\R^{m\times n}$, $b\in\R^{p-m}$ and $d\in\R^{m}$. Set
$$f(x)=\frac{1}{2}\|Ax-b\|^2  ~\text{and~}~h(x)=Cx-d.$$
The Lagrangian function of \eqref{CS12} is $$L(x,y)=\dfrac{1}{2}\|Ax-b\|^2-y^\top(Cx-d),$$ for any $x\in \S$ and $y\in\R^m$.
Direct calculations lead to
\begin{equation}\label{gra-hess-1}
\left\{
  \begin{array}{ll}
    \nabla h(x)=C,~\nabla_x L(x,y)=A^\top(Ax-b)-C^\top y,  & \hbox{ } \\
    \nabla^2f(x) = \nabla^2_{xx} L(x,y) = A^\top A, ~ \nabla^2 h(x) = 0,~\nabla^2 L(x,y) = \left[\begin{matrix}
A^\top A    & -C^\top\\
-C        &0
\end{matrix}\right]. & \hbox{ }
  \end{array}
\right.
\end{equation}

Since $\nabla^2 f(\cdot)$ and $\nabla^2 h(\cdot)$ are constant, Assumption \ref{assumption4} holds automatically everywhere. To ensure Assumptions 1$'$ and \ref{assumption3} hold, we introduce the following assumption on the input matrices $A$ and $C$.

\begin{assumption}\label{ACS}
For any index set $T\in \J_s$, $A_T$ is full column rank and $C_T$ is full row rank.
\end{assumption}
Suppose that Assumption \ref{ACS} holds.
Note that for any $(x,y)\in\R^{n+m}$ and $\beta>0$, $\T(x,y;\beta)\subseteq \J_s$. Together with rank$(\nabla_T h(x)) = \text{rank}(C_T)$, we can conclude that Assumption 1$'$ holds everywhere once $C_T$ is full row rank for all $T\in \J_s$. Similarly, since $\left(\nabla^2_{xx} L(x,y)\right)_{T,T} = A_T^\top A_T$, it is positive definite in the entire space $\R^s$ once $A_T$ is full column rank. Thus, Assumption \ref{assumption3} follows.

It is worth mentioning that Assumption \ref{ACS} is actually a mild condition for problem \eqref{CS12}. Indeed, the full column rankness of $A_T$ is the so-called $s$-regularity introduced by Beck and Eldar \cite{beck2013sparsity} which has been widely used in the CS community, and limiting the number of hard constraints will make the full row rankness of $C_T$ accessible (here $m$ is no more than $s$ and hence $s+m\leq 2s$). Under Assumption \ref{ACS}, we have the following optimality conditions for problem \eqref{CS12}.

\begin{proposition}\label{CS-optimality} Assume that the feasible set of problem \eqref{CS12} is nonempty and Assumption \ref{ACS} holds. Then the optimal solution set $S^*_{cs}$ of \eqref{CS12} is nonempty. Furthermore, for any strong $\beta$-Lagrangian stationary point $x^*$ of \eqref{CS12}, it is either a strictly local minimizer if $\|x^*\|_0 = s$ or a global optimal solution otherwise.
\end{proposition}
\noindent{\bf Proof.} The nonemptiness of $S^*_{cs}$ follows from the Frank-Wolfe Theorem and the observation $\S = \cup_{J\in \J_{s}}\R^n_J$. Then the rest of the assertion follows from Theorem \ref{sufficient-opt} and \cite[Theorem 4.2]{pan2017optimality}. \qed

\vspace{2mm}

With the above optimality results, we can apply LNA to solve \eqref{CS12} efficiently, since the linear system in each iteration is of size no more than $2s \times 2s$ and the algorithm will have a fast quadratic convergence rate, as stated in the following proposition.

\begin{proposition} \label{Cor-CS}
Suppose that Assumption \ref{ACS} holds. For given $\beta>0$, let $x^*$ be a strong $\beta$-Lagrangian stationary point of \eqref{CS12}	with $y^*$. Suppose that the initial point $z^0$ of sequence $\{z^k\}$ generated by LNA satisfies $z^0\in \N_S(z^*,\delta)$
where $\delta = \min\{\tilde{\delta}^*,1\}$ and $\tilde{\delta}^*$ is defined as in \eqref{delta-M}. Then for any $k\geq 0$,
\begin{itemize}
\item [(i)]$\lim\limits_{k\rightarrow \infty} z^k=z^*$ with quadratic convergence rate, i.e., $\|z^{k+1}-z^*\|\leq \dfrac{1}{2} \| z^k-z^* \|^2$.
\item [(ii)] LNA terminates with accuracy $\epsilon$ when $k\geq \left\lceil \frac{\log_2\left(4\delta^2 \sqrt{\left\|[A^\top A, -C^\top]\right\|^2+1}/\epsilon\right)}{2}\right\rceil$. 
\end{itemize}
\end{proposition}
\noindent{\bf Proof.} Since $\nabla^2 L(x,y)$ is constant and hence \eqref{A5} holds everywhere for any $L_2>0$. Thus, take $L_2 = 1/M^*$ with $M^*$ defined as in \eqref{delta-M}. Similarly, from \eqref{gra-hess-1}, we also have \eqref{A2} at any $z\in \R^{n+m}$ with $L_1 = \|[A^\top A, -C^\top]\|$.
By employing Theorem \ref{Thm-Qudratic-Convergence}, we can obtain the desired assertions.\qed

\subsection{Sparse High-Order Portfolio Selection}
In real financial markets, returns have often been found to be skewed and extreme events observed to be frequent which can be measured by skewness and kurtosis. Thus, based on Markowitz's mean-variance (MV) portfolio model, several studies consider the high-order portfolio selection with only a limited number of assets, i.e., the mean-variance-skewness-kurtosis model with cardinality constraint (MVSKC).
Suppose that $\tilde{r}\in\R^n$ is the return vector of $n$ assets and $x\in\R^n$ is the vector of portfolio weights.
The MVSKC model, which is first introduced in \cite{Boudt2015Higher}, takes the form of
\begin{equation}\label{MVSKC}
\begin{array}{ll}
\min ~&-\lambda_1x^\top\mu +\lambda_2 x^\top\Sigma x-\lambda_3x^\top\Phi (x\otimes x)+\lambda_4x^\top\Psi (x\otimes x\otimes x),\\
\text{s.t.}~&e^\top x=1,~\|x\|_0\leq s.
  \end{array}
\end{equation}
Here $\mu=E(\tilde{r})$ is the mean return vector, $\Sigma=E(rr^\top)$ is the covariance matrix, $\Phi=E(r(r^\top \otimes r^\top))$ is the co-skewness, $\Psi=E(r(r^\top \otimes r^\top \otimes r^\top))$ is the co-kurtosis matrix, with $r:=\tilde{r}-\mu$  the centered return vector and $\otimes$ the Kronecker product, $\lambda_i>0$, $i=1,\cdots,4$, are parameters to balance the four moments of the portfolio return. 
Set
$$f(x) = -\lambda_1x^\top\mu +\lambda_2 x^\top\Sigma x-\lambda_3x^\top\Phi (x\otimes x)+\lambda_4x^\top\Psi (x\otimes x\otimes x),~h(x) = e^\top x-1.$$
It is obvious that the objective function $f$ is nonconvex and twice continuously differentiable.
The corresponding Lagrangian function associated with problem \eqref{MVSKC} is
$$L(x,y) = f(x)-y(e^\top x-1), ~~~\forall x\in \S,~~y\in \R.$$
Utilizing Lemma 1 in \cite{wang2020sparse}, we have
\begin{equation}\label{gra-hess-2}
\left\{
  \begin{array}{ll}
    \nabla h(x)= e^\top,~ \nabla^2 h(x) = 0\\
    \nabla_x L(x,y)= -\lambda_1\mu +2\lambda_2 \Sigma x-3\lambda_3 \Phi (x\otimes x)+4\lambda_4\Psi (x\otimes x\otimes x)-y e,  & \hbox{ } \\
    \nabla^2f(x) = \nabla^2_{xx} L(x,y) = 2\lambda_2 \Sigma-6\lambda_3 \Phi (I\otimes x)+12\lambda_4\Psi (I\otimes x\otimes x),& \hbox{ } \\
    \nabla^2 L(x,y) = \left[\begin{matrix}
\nabla^2f(x)          &~~ -e \\
-e^\top    &~~0
\end{matrix}\right]. & \hbox{ }
  \end{array}
\right.
\end{equation}

It is easy to verify that Assumption 1$'$ holds directly for any $x\in\R^n$. Since $\nabla^2 f(\cdot)$ is continuously differentiable and $\nabla^2 h(\cdot)$ is constant, Assumption \ref{assumption4} automatically holds near $x^*$. To make the required Assumption \ref{assumption3} hold for problem \eqref{MVSKC}, we introduce the following assumption.

\vspace{2mm}

\begin{assumption}\label{S-portfolio}
$\lambda_i>0$, $i=1,\cdots,4$, satisfy $4\lambda_4(2\lambda_2-1)>\lambda_3^2$, 
and $\Sigma$ is positive definite restricted to the set $\{d\in\R^n:e^\top d=0\}$.
\end{assumption}
Learning from \eqref{gra-hess-2}, for any $d\in\{d\in\R^n:e^\top d=0\}$, we have
\begin{eqnarray*}
d^\top \nabla^2_{xx} L(x,y) d&=& d^\top (2\lambda_2 \Sigma-6\lambda_3 \Phi (I\otimes x)+12\lambda_4\Psi (I\otimes x\otimes x))d\\
~&=&d^\top E((2\lambda_2 -6\lambda_3r^\top x +12\lambda_4(r^\top x)^2 ))rr^\top)d\\
~&\geq & d^TE(rr^\top)d>0,
\end{eqnarray*}
where the second equality is from the definitions of $\Sigma$, $\Phi$ and $\Psi$, and the first inequality is from Assumption \ref{S-portfolio}.
Thus, Assumption \ref{assumption3} is valid.

Note that the condition in Assumption \ref{S-portfolio} is mild in real-world instances of sparse portfolio, since the covariance matrix is always positive definite. Under Assumption \ref{S-portfolio}, we have the following optimality conditions for problem \eqref{MVSKC} by Theorem \ref{optimalsta} and \cite[Theorem 4.2]{pan2017optimality}.
\begin{proposition}\label{portfolio-optimality}

(i) Suppose that $x^*$ is a local minimizer of  \eqref{S-portfolio}, then there exists a unique $y^*\in\R^m$ such that for any $\beta\in (0, \hat{\beta})$, $x^*$ is a strong $\beta$-Lagrangian stationary point, where $\hat{\beta}$ is defined as in \eqref{hatbeta}.

(ii) Assume that Assumption \ref{S-portfolio} holds and $x^*$ is a strong $\beta$-Lagrangian stationary point of \eqref{MVSKC}, then $x^*$ is a strictly local minimizer of \eqref{MVSKC}.
\end{proposition}

With the above optimality results, we can apply LNA to solve \eqref{MVSKC} efficiently. It is noteworthy that the computational cost per iteration is super low since we just need to handle $s+1$ linear equations in each iteration, and the algorithm will have a fast quadratic convergence rate  as stated below.

\begin{proposition}\label{Portfolio-Lemma}

Suppose that Assumption \ref{S-portfolio} holds.  For given $\beta>0$, let $x^*$ be a strong $\beta$-Lagrangian stationary point of \eqref{MVSKC} with $y^*$. Then LNA for \eqref{MVSKC} has locally quadratic convergence rate as stated in Theorem \ref{Thm-Qudratic-Convergence}.
\end{proposition}

\section{Numerical Experiments}
This section reports numerical results of  LNA in compressed sensing problem and sparse high-order portfolio selection on both synthetic and real data. All experiments were conducted by using MATLAB (R2018a)  on a laptop of 8GB memory and Inter(R) Core(TM) i5 1.8Ghz CPU.
 We terminate our method at $k$th step if $\eta_{\beta}(x^k,y^k;T_k) \leq 10^{-6}$ where $\eta_{\beta}(x^k,y^k;T)$ is defined as \eqref{Tolerance-Function} or $k$ reaches 1000.

\subsection{Compressed Sensing}
The aim of this subsection is to compare LNA with six state-of-the-art methods for compressed sensing problem \eqref{CS12}, including HTP\cite{foucart2011hard}\setcounter{footnote}{0}\footnote{
HTP is available at: \url{https://github.com/foucart/HTP.}}, NIHT\cite{blumensath2010normalized}\footnote{NIHT, GP and OMP are available at  \url{https://www.southampton.ac.uk/engineering/ about/staff/tb1m08.page$\#$software}. We use the version sparsify\_0\_5 in which NIHT, GP and OMP are called hard\_l0\_Mterm, greed\_gp and greed\_omp.\label{2}}, GP\cite{blumensath2008gradient}\footref{2}, OMP \cite{pati1993orthogonal,tropp2007signal}\footref{2}, CoSaMP \cite{needell2009cosamp}\footref{3} and  SP\cite{dai2009subspace}\footnote{
CoSaMP and  SP are available at: \url{http://media.aau.dk/null\_space\_pursuits/2011/07/a-few-corrections-to-cosamp-and-sp-matlab.html.}\label{3}}.

\subsubsection{Testing examples}
We generate the sensing matrix $\mathcal{A}$
in the same way as \cite{yin2015minimization,zhou2016null}.
Each column of $\mathcal{A}$ is normalized to $\|\mathcal{A}_j\|=1$ to make it consistent with the  algorithms used in \cite{foucart2011hard,blumensath2010normalized,blumensath2008gradient}. The true signal $x^*$ and the measurement  $\mathcal{B}$ are produced by the following pseudo MATLAB codes:
\begin{eqnarray*}
x^{*} = \verb"zeros"(n,1),~\Gamma = \verb"randperm"(n),~
x^{*}(\Gamma(1 : s)) = \verb"randn"(s, 1),~\mathcal{B}= \mathcal{A} x^*.
\end{eqnarray*}
Then, we randomly choose $m=\lceil 0.1s\rceil$ rows of $\mathcal{A}$ as $C$ in \eqref{CS12}. The rest part of $\mathcal{A}$ composes $A$ in the objective function. See the following pseudo MATLAB code for details:
\begin{eqnarray*}
&&J = \verb"randperm"(p),~J_1 = J(1:m),~J_2 = J(m+1:end);\\
&&A=\mathcal{A}(J_1),~b=\mathcal{B}(J_1),~C=\mathcal{A}(J_2),~d=\mathcal{B}(J_2).
\end{eqnarray*}

\begin{example}[Gaussian matrix]\label{ex-gau}
Let $\mathcal{A}\in \R^{p\times n}$  be a random Gaussian matrix  with each column
being identically and independently generated from the standard normal distribution.
\end{example}
\begin{example}[Partial DCT  matrix]\label{ex-dct}
Let $\mathcal{A}\in \R^{p\times n}$  be a random partial discrete cosine transform (DCT) matrix generated by
\begin{eqnarray*}
\mathcal{A}_{ij}= \cos(2\pi (j-1) \psi_i),~~i=1,\ldots, p,~~j=1,\ldots, n
\end{eqnarray*}
where $\psi_i$ ($i=1,\ldots, m$) is uniformly and independently sampled from $[0,1]$.
\end{example}

\subsubsection{Numerical comparisons}
We set the maximum number of iterations and the tolerance as $1000$ and $10^{-6}$, respectively, in all of the six comparison methods mentioned above. The initializations are set to be $x^0=0$, $y^0=0$ and $\beta=5/n$ for LNA.
For comparison purpose, HTP, NIHT, GP, OMP, CoSaMP and SP as tested in this section are all initialized with the origin in their default setups.

We say a recovery of a method is successful if $\|x-x^*\|<0.01\|x^*\|$, where $x$ is the solution produced by this method. The corresponding success rate is defined as the percentage of the number of successful recovery instances over all trials.

\begin{figure}
	\centering  
	\subfigure{
		\includegraphics[width=0.42\linewidth]{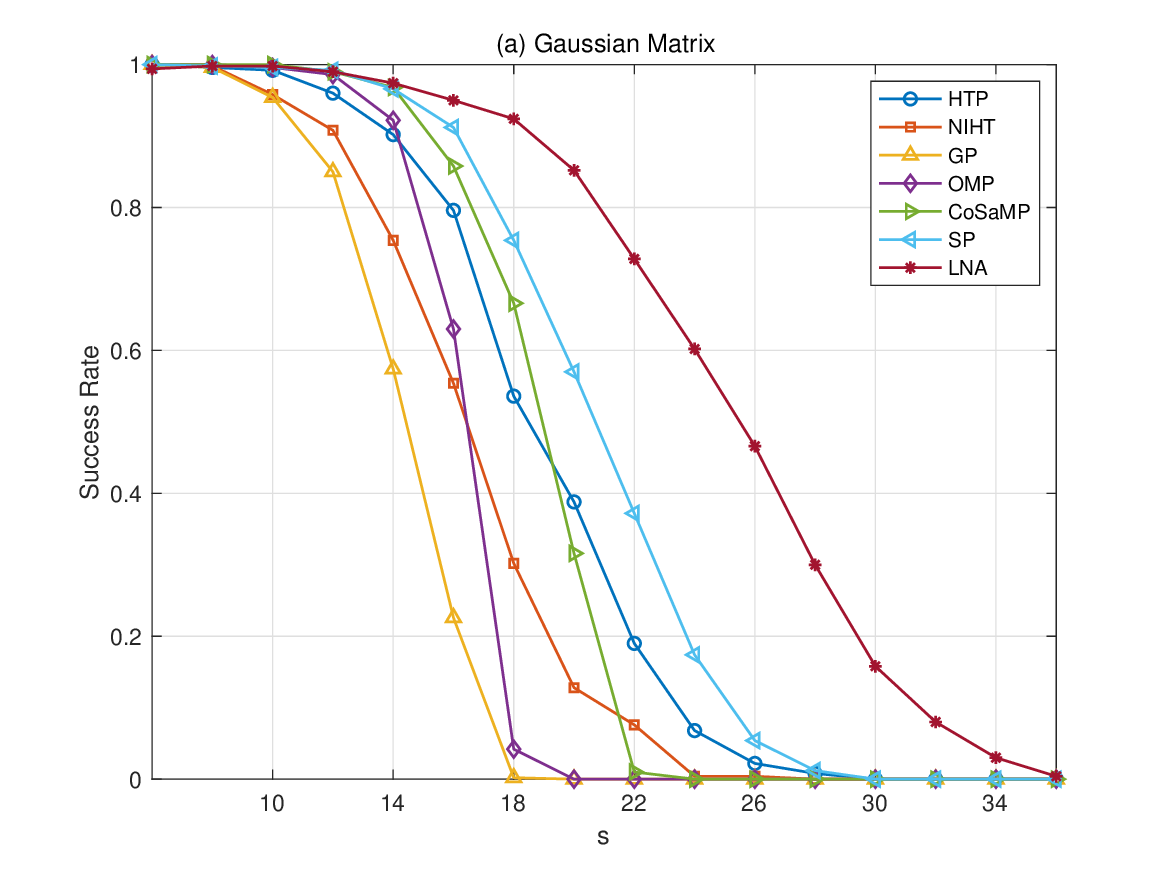}}
	\subfigure{
		\includegraphics[width=0.42\linewidth]{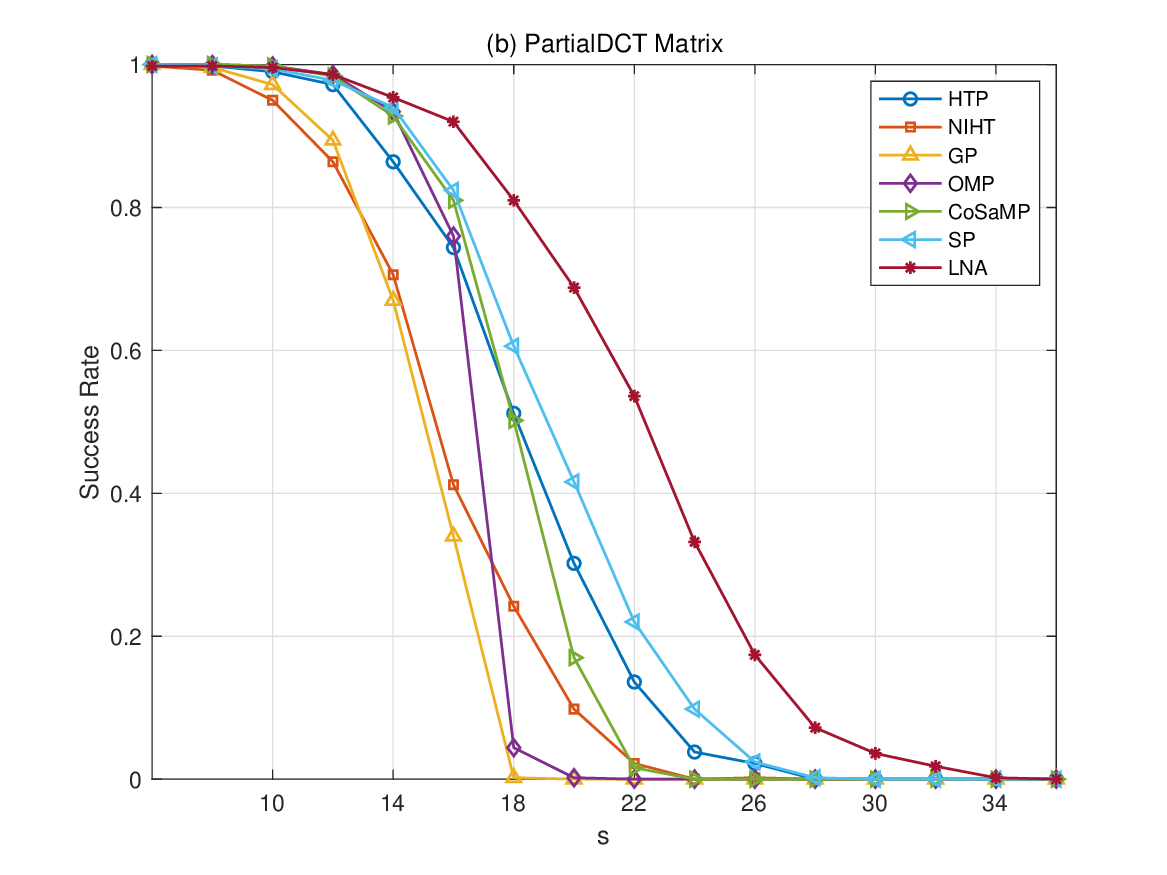}}
		\caption{Success rates. $n = 256, p=\lceil n/4\rceil,  s\in\{6,8,\ldots, 36\}$.}
\label{fig:SuccRate}
	\subfigure{
		\includegraphics[width=0.42\linewidth]{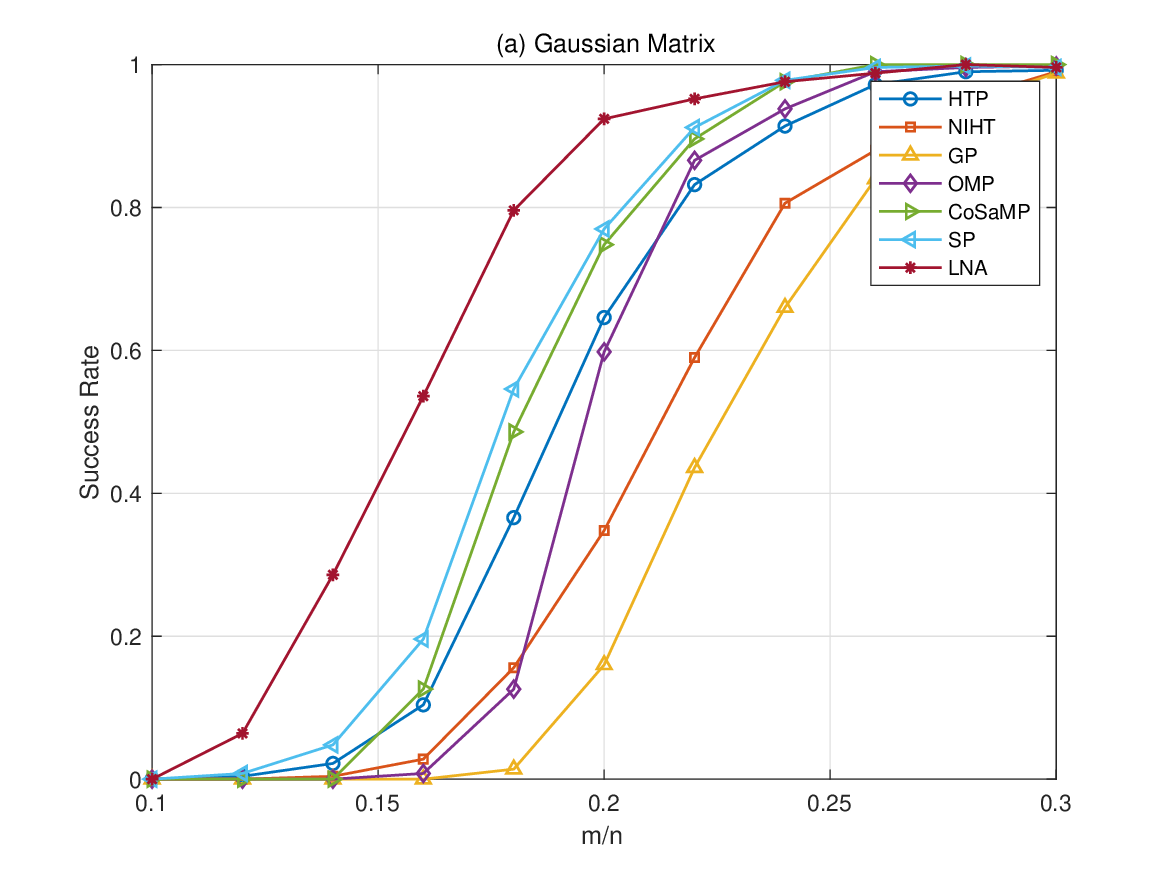}}
	\subfigure{
		\includegraphics[width=0.42\linewidth]{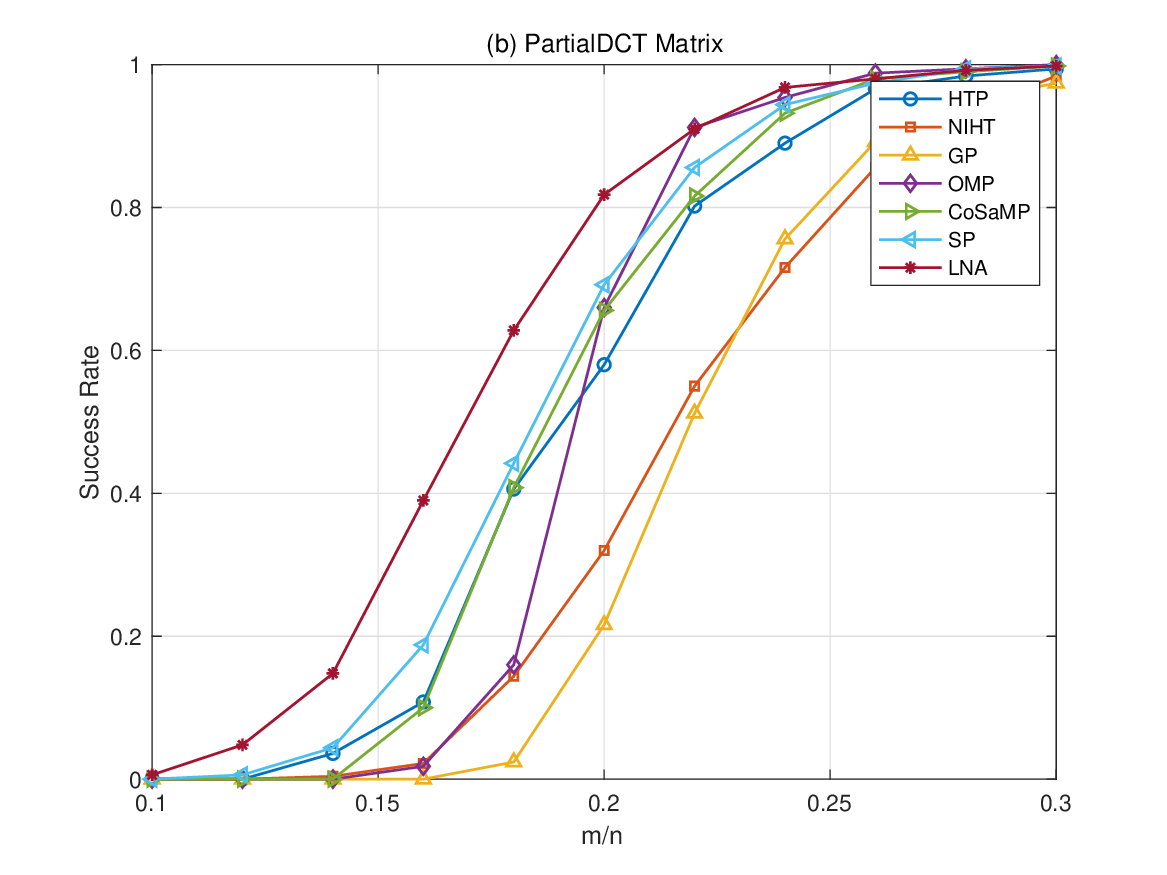}}
		\caption{Success rates. $n = 256, s=\lceil 0.05 n\rceil, p=\lceil rn\rceil$ with $r\in\{  0.1, 0.12, \ldots, 0.3\}$.}
\label{fig:SuccRate-mn}
\end{figure}

Firstly, we run $500$ independent trials with fixed $n=256, p=\lceil n/4\rceil$ at different sparsity levels $s$ from 6 to 36. The corresponding success rates are illustrated in Fig. \ref{fig:SuccRate}. One can observe that LNA always yielded the highest success rate for each $s$ under both Example \ref{ex-gau} and Example \ref{ex-dct}, while a lowest success rate is generated in {GP}.
For example, when $s=20$ for Gaussian matrix, $85\%$ successful recoveries are guaranteed in our method, which performed much better than other methods, whose success rates are all less than $60\%$.

Next, we implement $500$ independent trials by varying $p=\lceil r n\rceil$ in $r\in\{0.1,0.12,\ldots,0.3\}$ when $n=256, s=\lceil 0.05 n\rceil$ in Fig. \ref{fig:SuccRate-mn}, which indicates that the larger $m$ is, the easier the problem becomes to be solved.
Again, LNA outperformed the others in the success rate for each $s$, and {GP} still came the last.

We now examine these algorithms with higher dimensions $n$ between 5000 and 25000 with 50 trials when $p=\lceil n/4\rceil, s=\lceil 0.01n\rceil, \lceil0.05n\rceil$ in the framework of Example \ref{ex-gau}, to compare their speed of convergence and the accuracy of solutions. The average absolute error $\|x-x^*\|$ and CPU time are  presented in Table \ref{tab:error} and Table \ref{tab:time}, respectively. One can see that the highest accurate recovery can be obtained in LNA with the least CPU time for most cases. Although OMP and  HTP  rendered solutions as accurate as those by LNA when $n$ is small, they presented some shortcomings by comparing to LNA. In OMP, the accuracy cannot be guaranteed when $n$ is large, and some inaccurate ones were produced when $s=\lceil0.05n\rceil$ and $n\geq 20000$ particularly, which implies that OMP only worked well when the solution is very sparse. On the other hand, the CPU time consumed by HTP is booming over $n$. For example, when $n=25000$ and $s=\lceil0.05n\rceil$,  5.99 seconds  by LNA against 159.38 seconds by HTP. Moreover, even though NIHT is the fastest one among the six  methods, its accuracy is much worse than others as it is stable at achieving the solutions with accuracy of order $10^{-7}$. That is to say, the superiority of LNA becomes more obvious in the trade off of high accuracy and convergence speed with high dimensional data.

 \begin{table}\centering
 \caption{Average absolute error $\|x-x^*\|$ for Example   \ref{ex-gau}.
 \label{tab:error}}
\begin{tabular}{lllllllll}
\hline\noalign{\smallskip}
$s$&$n$&{LNA} &  {HTP}&{NIHT}&{GP}&{OMP}&{CoSaMP} & SP \\
\noalign{\smallskip}\hline\noalign{\smallskip}
\multirow{5}{*}{$\lceil0.01n\rceil$} &	5000	&	2.71e-15	&	3.13e-15	&	2.03e-8	&	4.01e-15	&	2.78e-15	&	1.41e-14	&	1.41e-14	\\
	&	10000	&	4.86e-15	&	5.70e-15	&	2.27e-8	&	7.04e-15	&	4.80e-15	&	2.15e-14	&	2.15e-14	\\
	&	15000	& 6.52e-15	&	7.39e-15	&	2.92e-8	&	1.06e-14	&	6.82e-15	&	2.98e-14	&	2.98e-14	\\
	&	20000	&	8.87e-15	&	9.97e-15	&	4.37e-8	&	1.37e-14 &	9.34e-15	&	4.08e-14	&	4.08e-14	\\
	&	25000	&	1.04e-14	&	1.21e-14	&	3.95e-8	&	1.72e-14	&	1.16e-14	&	4.44e-14	&	4.44e-14	\\
	\noalign{\smallskip}\hline\noalign{\smallskip}
\multirow{5}{*}{$\lceil0.05n\rceil$}	&	5000	&	1.14e-14	&	1.11e-14	&	1.63e-7	&	1.49e-14	&	1.08e-14&	4.10e-14&	4.10e-14	\\
	&	10000	&	2.80e-14	&	2.28e-14	& 3.30e-7	&	2.97e-14	&	2.42e-14&	7.75e-14&	7.75e-14	\\
	&	15000	&	3.91e-14	&	3.70e-14	&	3.06e-7	&	5.02e-14	&	4.31e-14&	1.10e-13&	1.10e-13	\\
	&	20000	&	5.22e-14	&	4.77e-14	&	4.03e-7	&	5.83e-14	&	5.15e-04&	1.34e-13&	1.34e-13	\\
	&	25000	&	6.30e-14	&	6.12e-14	&	3.75e-7	&	7.74e-14	&	6.30e-04&	1.82e-13&	1.82e-13	\\
	\noalign{\smallskip}\hline
    \end{tabular}
    \vspace{3mm}
 \caption{Average CPU time (in seconds) for Example  \ref{ex-gau}.
  \label{tab:time}}
\begin{tabular}{lllllllll}
\hline\noalign{\smallskip}
$s$&$n$&{LNA} &  {HTP}&{NIHT}&{GP}&{OMP}&{CoSaMP} & SP \\
\noalign{\smallskip}\hline\noalign{\smallskip}
\multirow{5}{*}{$\lceil0.01n\rceil$} &	5000	&	0.06	&	0.62&	0.21	&	1.92&	0.30&	0.57&	0.04	\\
	&	10000	&	0.25	&	4.02	&	0.83	&	14.07	&	2.28	&	0.28	&	0.18	\\
	&	15000	&	0.60	&	13.07	&	1.93	&	46.88&	7.70	&	1.01	&	0.78	\\
	&	20000	&	1.08	&	32.08	&	3.50	&	110.07	&	18.16	&	2.13&	1.34	\\
	&	25000	&	1.77	&	111.94	&	6.11	&	230.27&	37.43&	4.11	&	2.57\\
	\noalign{\smallskip}\hline\noalign{\smallskip}
\multirow{5}{*}{$\lceil0.05n\rceil$}	&	5000	&	0.14	&	0.80&	0.62	&	2.18	&	1.73	&	1.44	&	0.93\\
	&	10000	&	0.59	&	5.71	&	2.42&	15.07	&	13.66	&	13.87	&	5.36	\\
	&	15000	&	1.50	&	18.90	&	5.47	&	50.41	&	46.82	& 47.41	&	18.23	\\
	&	20000	&	3.03	&	49.10	&	11.24	&	118.35	&	111.40	&	127.37	&	55.21	\\
	&	25000	&	5.99	&	159.38&	17.49	&	239.70	&	174.21	&	217.07	&	96.38\\
	\noalign{\smallskip}\hline
    \end{tabular}
\end{table}

\begin{figure}
	\centering
		\includegraphics[scale=0.33]{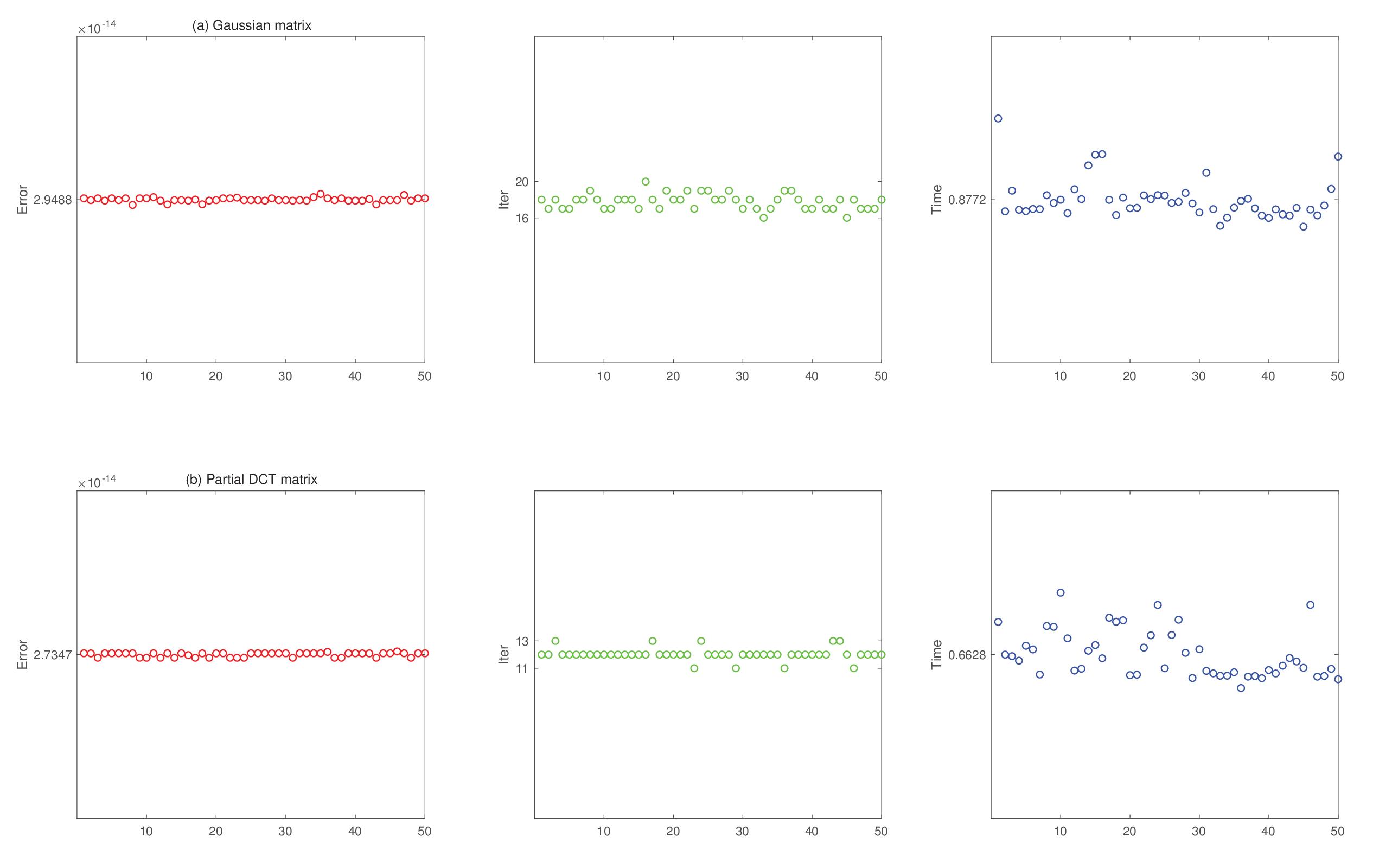}
		\caption{Effects of initial points for Examples \ref{ex-gau} and \ref{ex-dct}.}
\label{fig:Robust x0}
\end{figure}

As stated in Theorem \ref{Thm-Qudratic-Convergence}, LNA is a local method. Therefore,
We conduct numerical experiments with randomly generated initial points for CS problems to see how the initial points would affect LNA.
To proceed, we apply LNA into solving Examples \ref{ex-gau} and \ref{ex-dct} with $n=10000, p=\lceil n/4\rceil, s=\lceil 0.05 n\rceil$.
We run the LNA under 50 different initial points which are randomly generated from the uniform distribution, namely, $(x^0,y^0)\thicksim\text{U}[0,1]$.
The absolute error $\|x-x^*\|$, the number of iterations and CPU time are plotted in Fig. \ref{fig:Robust x0}, where the x-axis stands for the 50 initial points. One can see that all the results stabilize at a certain level, which indicates that LNA is not sensitive to the choices of the initial points for CS problems.

\subsection{Sparse High-Order Portfolio Selection}

This subsection is devoted to comparing LNA with successive convex approximation algorithm (SCA)\cite{wang2020sparse} in sparse high-order portfolio selection problem \eqref{MVSKC} on real data sets.
\subsubsection{Testing examples}

\begin{example}\label{portfolio-example} (Portfolio data sets)
The data sets used in our experiments are selected from the Standard and Poor's 500 (USA) (S\&P 500 for short)\setcounter{footnote}{5}\footnote{ \url{http://cran.r-project.org/web/packages/portfolioBacktest/vignettes }\label{6}}. Firstly, we randomly select 100 socks from S\&P 500 Index components and randomly choose 500 continuous trading days from  2012-12-01 to 2018-12-01. Then the selected data is normalized to raise precision of the model, and the sample moments\setcounter{footnote}{6}\footnote{ \url{http://www.mathworks.com/matlabcentral/fileexchange/47839-co\_moments-m.}\label{7}} are computed. 
To be fair enough, we adopt the choices of model parameters in \eqref{MVSKC} from \cite{Boudt2015Higher} with $\lambda_1=1,~\lambda_2=\xi/2,~\lambda_3=\xi(\xi+1)/6,~\lambda_4=\xi(\xi+1)(\xi+2)/24$, where $\xi=5,~10$ is the risk aversion parameter. Direct calculations certify that these parameters satisfy the condition in Assumption \ref{S-portfolio}. Additionally, the sparsity level $s$ for $\xi=5$, $10$ will be varying among $\{5, 10, 15, 20, 25\}$ to generate a total of $10$ testing instances.
\end{example}

\subsubsection{Numerical comparisons}

For portfolio data sets in Example \ref{portfolio-example}, we found that different initial points lead to different output solutions. This is reasonable since the objective function is nonconvex and LNA is a locally convergence method. Note that there are various ways to find an initial point near to strong $\beta$-Lagrangian stationary point, for instance, some first-order gradient descent methods and convex relaxation methods.
For simplicity, we initialize LNA with the origin $(x^0, y^0)=(0,0)$ which is appropriate to our testing examples, and $\beta=1$. For comparison purpose, SCA is called with the initial point $x^0=0$, along with other parameters as $\alpha=0.2$, $\rho=4\times 10^{-3},3\times 10^{-3}$ for $\xi=5,10$ respectively. The sparsity of a solution $x$ generated by SCA will be recorded by $\hat{s}:=\min\{t:\sum_{i=1}^t |x|_{(i)}\geq 0.99\|x\|_1\}$.
Table \ref{tab:objtime} records $s$, $\hat s$, the objective function value (f-value) and CPU time when $\xi=5$, $10$ respectively.

\begin{table}[h]
\centering
 \caption{$\hat s$, f-value and CPU time (in seconds) for Example \ref{portfolio-example}.
 \label{tab:objtime}}
\begin{tabular}{l|l|ll|ll|l|ll|ll}
\hline
                     & \multicolumn{5}{c|}{$\xi=5$}                                                                                   & \multicolumn{5}{c}{$\xi=10$}                                                                                  \\ \hline
\multirow{2}{*}{$s$} & \multicolumn{1}{c|}{\multirow{2}{*}{$\hat{s}$}} & \multicolumn{2}{c|}{f-value} & \multicolumn{2}{c|}{CPU time} & \multicolumn{1}{c|}{\multirow{2}{*}{$\hat{s}$}} & \multicolumn{2}{c|}{f-value} & \multicolumn{2}{c}{CPU time} \\ \cline{3-6} \cline{8-11}
                     & \multicolumn{1}{c|}{}                           & LNA           & SCA          & LNA           & SCA           & \multicolumn{1}{c|}{}                           & LNA           & SCA          & LNA           & SCA           \\ \hline
5                    & 55                                              &              -1.11 &      -1.50        &      15.55         &       50.55        & 88                                              &     -0.70          &             -1.44 &     95.73          &      144.92         \\
10                   & 55                                              &              -3.08 &        -1.50        &       39.53          &     45.48      & 88                                              &      -1.55         &-1.44              &    46.66           &     216.47          \\
15                   & 55                                              &              -2.81 &     -1.50         &      26.80         &      45.43         & 88                                              &     -1.78          &             -1.44 &       77.15        &       186.79        \\
20                   & 55                                              &              -4.27 &      -1.50        &       23.74        &       47.06        & 88                                              &     -2.80          &             -1.44 &       75.58        &        194.80       \\
25                   & 55                                              &              -4.35 &     -1.50         &       34.30        &      46.62         & 88                                              &     -3.01          &              -1.44     &    23.96      &         117.98      \\ \hline
\end{tabular}
\end{table}

As one can see from Table \ref{tab:objtime}, LNA outperforms SCA in computational time for all testing instances, and attains smaller f-value than that of SCA when $s>5$. Specifically, LNA provides $s$-sparse solutions while SCA fails. 
Furthermore, as $s$ grows, f-value decreases in LNA, which indicates the trend to the true minimum in some sense, against that almost no improvement of f-value in SCA.

\section{Conclusion}
In this paper, we have designed a second-order greedy algorithm named the Lagrange-Newton Algorithm (LNA) for the sparse nonlinear programming (SNP) problem with sparsity and nonlinear equality constraints, based on the
strong $\beta$-Lagrangian stationarity and Lagrangian equations.
The resulting LNA has shown to be effective, with local quadratic convergence rate and low iterative complexity from the theoretical perspective, and good computational superiority  from the numerical perspective.

There are also some issues that remain to be further investigated. As LNA is a second-order local method with heavy reliance on the initial points in general, the first issue is whether a line search scheme would be equipped for LNA, attempting to achieve global convergence and to accelerate the algorithm. Another more general issue would be whether we can extend LNA to more general sparse optimization models with equality and inequality constraints. We leave these in our future research.

\section*{Acknowledgement}
We would like to thank AE and two referees for their valuable comments to improve our paper, and Dr. Shenglong Zhou for his great support on the numerical experiments.

\bibliographystyle{spmpsci}   
\bibliography{bib_LNA}
\end{document}